\documentclass[a4paper,11pt]{article}
\usepackage{amssymb}
\newtheorem{theorem}{Theorem}[section]
\newtheorem{lemma}{Lemma}[section]
\newtheorem{corollary}{Corollary}[section]
\newtheorem{definition}{Definition}[section]
\newtheorem{proposition}{Proposition}[section]

\newenvironment{remark}{\vspace{1ex}{\bf Remark.}\rm}{\vspace{1ex}}
\newenvironment{remarks}{\vspace{1ex}{\bf Remarks.}\rm}{\vspace{1ex}}
\newenvironment{proof}{{\bf Proof.}}{\par\hspace{25em}\rule{1ex}{1ex}\par}
\newcommand{\nablast}{\stackrel{*}{\nabla}}
\title{Transverse totally geodesic submanifolds of the tangent
bundle\footnote{Math. Publ. Debrecen 64/1-2 (2004), 129-154}.}
\author{Mohamed Tahar Kadaoui ABBASSI\\
Alexander YAMPOLSKY}
\date{July 14, 2002}
\begin{document}
\maketitle
\begin{abstract}
It is well-known that if $\xi$ is a smooth vector field on a given
Riemannian manifold $M^n$ then $\xi$ naturally defines a
submanifold $\xi(M^n)$ transverse to the fibers of the tangent
bundle $TM^n$ with Sasaki metric. In this paper, we are interested
in transverse totally geodesic submanifolds of the tangent bundle.
We show that a transverse submanifold $N^l$ of $TM^n$ ($1 \leq l
\leq n$) can be realized locally as the image of a submanifold
$F^l$ of $M^n$ under some vector field $\xi$ defined along $F^l$.
For such images $\xi(F^l)$, the conditions to be totally geodesic
are presented. We show that these conditions are not so rigid as
in the case of $l=n$, and we treat several special cases ($\xi$ of
constant length, $\xi$ normal to $F^l$, $M^n$ of constant
curvature, $M^n$ a Lie group and $\xi$ a left invariant vector
field).\\[2ex]
%[1ex]
{\it AMS subject class:} Primary 53B25, 53C07, 53C40; Secondary
22E15, 53B20, 53B21, 53C12, 53C25, 53M10.

\end{abstract}
\section*{Introduction.}

Let $(M^n,g)$ be a Riemannian manifold and $(TM^n,g_s)$ its
tangent bundle equipped with the Sasaki metric \cite{Sk}. Let
$\xi$ be a given smooth vector field on $M^n$. Then $\xi$
naturally defines a mapping $\xi:M^n\to TM^n$ such that the
submanifold $\xi(M^n)\subset TM^n$ is transverse to the fibers.
This fact allows to ascribe to the vector field $\xi$ some
geometrical characteristics from the geometry of submanifolds. We
say that the vector field $\xi$ is  \textit{minimal, totally
umbilic} or \textit{totally geodesic} if $\xi(M^n)$ possesses the
same property. In a similar way we can say about the
\textit{sectional, Ricci} or \textit{scalar curvature} of a vector
field. For the case of a \textit{unit} vector field this approach
has been proposed by H.Gluck and W.Ziller \cite{G-Z}. They proved
that the Hopf vector field $h$ on three-sphere $S^3$ is one with
globally minimal volume, i.e. $h(S^3)$ is a globally minimal
submanifold in the unit tangent bundle $T_1S^3$. Corresponding
local consideration leads to the notion of the \textit{ mean
curvature of a unit vector field} and a number of examples of
locally minimal unit vector fields were found based on a preprint
version of \cite{GM} (see \cite{BX-V1,BX-V2,GD-V1} and
references). In a different way, the second author found examples
of unit vector fields of \textit{constant mean curvature}
\cite{Ym1} and completely described the \textit{totally geodesic }
unit vector fields on 2-dimensional manifolds of constant
curvature \cite{Ym2}. The energy of a mapping $\xi:M^n\to T_1M^n$
can also be ascribed to the vector field $\xi$ and we can say
about the \textit{energy } of a unit vector field (see \cite{Wd,
GMLF, Wk-Bt} and references).

In contrast to unit vector fields, there are few results (both of
local or global aspects) on the geometry of general vector fields
treated as submanifolds in the \textit{tangent bundle}. It is
known \cite{Liu} that if $\xi$ is the zero vector field, then
$\xi(M^n)$ is totally geodesic in $TM^n$. Walczak~P. \cite{Wk}
treated the case when $\xi$ is a non-zero vector field on $M^n$
and proved that if $\xi$ is a parallel vector field on $M^n$, then
$\xi(M^n)$ is totally geodesic in $TM^n$. Moreover, if $\xi$ is of
constant length, then $\xi(M^n)$ is totally geodesic in $TM^n$ if
and only if $\xi$ is a parallel vector field on $M^n$. The latter
condition is rather burdensome. The basic manifold $M^n$ should be
a metrical product $M^{n-k}\times E^k\ (k\geq1)$, where $E^k$ is a
Euclidean (flat) factor.

Remark that $\xi(M^n)$ has maximal dimension among submanifolds in
the tangent bundle, transverse to the fibers. In this paper, we
study submanifolds $N^l$ of $TM^n$ with $l\le n$ which are
transverse to the fibers. We show in section 2 that any transverse
submanifold $N^l$ of $TM^n$ can be realized locally as the image
of a submanifold $F^l$ of $M^n$ under some vector field $\xi$
defined along $F^l$. We also investigate some cases when the image
can be globally realized. Mainly, we are interested in
submanifolds among this class which are totally geodesic. In this
way, we get a chain of inclusions:
$$\xi(F^l)\subset\xi(M^n)\subset TM^n. $$ In comparison with the
case when $\xi$ is defined over the whole $M^n$ or, at least, over
a domain $D^n\subset M^n$ as in \cite{Wk}, the picture becomes
different, because $\xi(F^l)$ can be totally geodesic in $TM^n$
while $\xi(M^n)$ is not. Our considerations include also the case
when the vector field is defined only on $F^l$, so that $\xi$
defines a ``direct" embedding $\xi:F^l\to TM^n$.

For $l=1$ we get nothing else but a vector field along a curve in
$M^n$ which generates a geodesic in $TM^n$. Sasaki S. \cite{Sk}
described geodesic lines in $TM^n$ in terms of vector fields along
curves in $M^n$ and found the differential equations on the curve
and the corresponding vector field. Moreover, in the case when
$M^n$ is of constant curvature, Sato K. \cite{St} explicitly
described the curves and the vector fields.

Evidently, our approach takes an intermediate position between the
above mentioned considerations for $l=1$ and $l=n$.

Necessary and sufficient conditions on $\xi(F^l)$ to be totally
geodesic, that we make explicit in section~3 (Proposition
\ref{Pr5}), have a clearer geometrical meaning if we suppose that
$\xi$ is of constant length along $F^l$ (Theorem \ref{Th1}) or is
a normal vector field along $F^l$ (Theorem \ref{Th2}). Indeed, an
application of Theorem \ref{Th2} to the specific case of foliated
Riemannian manifolds allows us to clarify the geometrical
structure of $\xi(M^n)$ (Corollary \ref{Foli}).

The case of a base space $M^n$ of constant curvature is discussed
in detail in section~4. An application to the case of a Riemannian
manifold of constant curvature enlightens us as to the non
rigidity of the totally geodesic property of $\xi(F^l)$, $l < n$,
contrary to the case $l=n$.

Finally, an application of our results to Lie groups endowed with
bi-invariant metrics gives a clear geometrical picture of our
problem.

\begin{remark} Throughout the paper
\begin{itemize}\itemsep=0ex
\item[-] $M^n$ is a given Riemannian manifold with metric $\bar g$,
$F^l$ is a submanifold of $M^n$ with the induced metric $g$,
$TM^n$ is the tangent bundle of $M^n$ equipped with the Sasaki
metric $g_s$;
\item[-] $\bar\nabla , \, \nabla , \, \tilde \nabla $ are the
Levi-Civita connections with respect to $\bar g,\, g,\, g_s$
respectively;
\item[-] the indices range is fixed as $a,b,c=1\dots n; \ \
i,j,k=1\dots l$;
\item[-] all the vector fields  are supposed sufficiently smooth,
say of class $C^\infty$.
\end{itemize}
\end{remark}

\section{Local geometry of $\xi (F^l)$.}

\subsection{Tangent bundle of $\xi(F^l)$.}

Let $(M^n,\bar g)$ be an $n$-dimensional Riemannian manifold  with metric $\bar g$.
Denote by $\bar g (\cdot\,,\cdot)$ the scalar product with respect to $\bar g$. The
{\it Sasaki metric} $g_s$ on $TM^n$ is defined by the following scalar product: if
$\tilde X,\tilde Y$ are tangent vector fields on $TM^n$,  then
\begin{equation}
\label{Eqn1}
   g_s(\tilde X,\tilde Y)=
   {\bar g}(\pi_* \tilde X, \pi_* \tilde Y)+{\bar g}(K \tilde X,K \tilde Y)
\end{equation}
where $\pi_*:TTM^n \to TM^n $ is the differential of the projection $\pi:TM^n \to M^n
$ and $K: TTM^n \to TM^n$  is the {\it connection map} \cite{Dmb}. The local
representations for $\pi_*$ and $K$ are the following ones. Let $(x^1,\dots ,x^n)$ be
a local coordinate system on $M^n$. Denote by $\partial/\partial x^a $ the natural
tangent coordinate frame. Then, at each point $x\in M^n$, any tangent vector $\xi$ can
be decomposed as $\xi=\xi^a \frac{\partial}{\partial x^a}(x)$. The set of parameters
$\{x^1,\dots ,x^n;\,\xi^1,\dots,\xi^n\}$ forms the natural induced coordinate system
in $TM^n$, i.e. for a point $z=(x,\xi )\in TM^n$, with $x\in M^n, \ \ \xi \in T_xM^n$,
we have $x=(x^1,\dots ,x^n), \, \xi =\xi ^a\frac{\partial}{\partial x^a}(x)$. The
natural frame in $T_{z}TM^n$ is formed by $\left\{ \frac{\partial}{\partial x^a}(z),
\frac{\partial}{\partial \xi ^a}(z)\right\}$ and for any $\tilde X\in T_{z}TM^n$ we
have the decomposition $\tilde X=\tilde X^a\frac{\partial}{\partial x^a}(z)+\tilde
X^{n+a}\frac{\partial}{\partial    \xi ^a}(z)$. Now locally, the \textit{horizontal}
and \textit{vertical} projections of $\tilde X$ are given by
   \begin{equation} \label{Eqn2}
   \begin{array}{l}
   \pi_* \tilde X= \tilde X^a\frac{\partial}{\partial x^a}(\pi(z)), \\[1ex]
   K \tilde X= (\tilde X^{n+a}+\bar\Gamma^a_{bc}(\pi(z))\,\xi^b \tilde X^c)\,
   \frac{\partial}{\partial x^a}(\pi(z)), \\[1ex]
   \end{array}
   \end{equation}
where $\bar\Gamma^a_{bc}$ are the Christoffel symbols of the metric $\bar g$. The
inverse operations are called \textit{lifts }. If $\bar X=\bar X^a\,\partial/\partial
x^a$ is a vector field on $M^n$ then the vector fields on $TM$ given by
   $$
   \begin{array}{l}
   \bar X^h=\bar X^a \partial/\partial x^a-\bar\Gamma^a_{bc}\,\xi^b\bar X^c\,\partial/\partial
   \xi^a ,\\[1ex]
   \bar X^v=\bar X^a\partial/\partial \xi^a
   \end{array}
   $$
are called the \textit{horizontal} and \textit{vertical} lifts of $X$ respectively.
Remark that for any vector field $\bar X$ on $M^n$ it holds
   \begin{equation}\label{Pr}
   \begin{array}{ll}
   \pi_* {\bar X}^h=\bar X,& K {\bar X}^h=0, \\[1ex]
   \pi_* {\bar X}^v=0,       & K {\bar X}^v=\bar X.
   \end{array}
   \end{equation}

Let $F^l$ be an $l$-dimensional  submanifold in $M^n$ with  a local representation
given by
$$
x^a= x^a(u^1,\dots ,u^l).
$$
Let $\xi$ be a vector field on $M^n$ defined in some neighborhood of (or only on) the
submanifold $F^l$. Then the restriction of $\xi$ to the submanifold $F^l$, called
\textit{a vector field on $M^n$ along $F^l$}, generates a submanifold $\xi(F^l)\subset
TM^n$ with a local representation of the form
\begin{equation}\label{Sub}
\xi(F^l):
\left\{
\begin{array}{ll}
x^a=& x^a(u^1,\dots ,u^l), \\
\xi ^a=&\xi ^a(x^1(u^1,\dots ,u^l),\dots ,x^n(u^1,\dots ,u^l)).
\end{array}
\right.
\end{equation}
In what follows we will refer to the submanifold (\ref{Sub}) as to one
\textit{generated by  a vector field on $M^n$ along $F^l$.}

The following Proposition describes the tangent space of $\xi(F^l)$.

    \begin{proposition}\label{Pr1}
    A vector field $\tilde X$ on $TM^n$ is tangent to $\xi(F^l)$ along $\xi(F^l)$
    if and only if its horizontal-vertical decomposition is of the form
    $$
    \tilde X = X^h+(\bar \nabla _X\, \xi)^v,
    $$
    where $X$ is a tangent vector field on $F^l$, $\bar \nabla _X\, \xi$ is the covariant
    derivative of $\xi$ in the direction of $X$ with respect to the Levi-Civita connection
    of $M^n$ and the lifts are considered as those on $TM^n$.
    \end{proposition}

\begin{proof}
Let us denote by $\tilde e_i$ the vectors of the coordinate frame
of $\xi (F^l)$. Then, evidently, $$\textstyle \tilde e_i=\left\{
\frac{\partial x^1}{\partial u^i}, \dots , \frac{\partial
x^n}{\partial u^i}; \ \ \frac{\partial \xi ^1}{\partial u^i},
\dots , \frac{\partial \xi ^n}{\partial u^i} \right\}. $$ Applying
(\ref{Eqn2}), we have $$
\begin{array}{rl}
\pi_*\tilde e_i= &\frac{\partial x^a}{\partial
u^i}\frac{\partial}{\partial x^a}= \frac{\partial}{\partial
u^i},\\[2ex] K\tilde e_i =& (\frac{\partial \xi^a}{\partial u^i} +
\bar \Gamma ^a_{bc} \, \xi^b \, \frac {\partial x^c}{\partial
u^i})\frac{\partial}{\partial x^a}= (\frac{\partial
\xi^a}{\partial x^c}\frac{\partial x^c}{\partial u^i} + \bar
\Gamma ^a_{bc} \, \xi^b \, \frac {\partial x^c}{\partial
u^i})\frac{\partial}{\partial x^a}\\[1ex] =&\frac {\partial
x^c}{\partial u^i}(\frac{\partial \xi^a}{\partial x^c} + \bar
\Gamma ^a_{bc} \, \xi^b \,) \frac{\partial}{\partial x^a}=
\bar\nabla_i\xi,
\end{array}
$$ where $\bar \Gamma ^a_{bc}$ are the Christoffel symbols of the
metric $\bar g$ taken along $F^l$ and  $\bar \nabla _i$ means the
covariant derivative of a vector field on $M^n$ with respect to
the Levi-Civita  connection of $\bar g$  along the $i$-th
coordinate curve of the submanifold $F^l \subset M^n$. Summing up,
we have
\begin{equation}\label{Eqn3}
\tilde e_i= \left(\frac{\partial}{\partial u^i}\right)^h +(\bar
\nabla _i \xi)^v.
\end{equation}

Let $\tilde X$ be a vector field on $TM^n$ tangent to $\xi(F^l)$
along $\xi(F^l)$. Then the following decomposition holds
$
\tilde X \,= \tilde X^i \tilde e_i.
$
Set $ X=\tilde X^i\partial/\partial u^i$. The vector field $X$ is
tangent to $F^l$ and, taking into account (\ref{Eqn3}), the
decomposition of $\tilde X$ can be represented as
$
\tilde X=X^h+(\bar \nabla _X \, \xi)^v,
$
which completes the proof.

\end{proof}

\begin{corollary}\label{Cor1} Let $(F^l,g)$ be a submanifold of a
Riemannian manifold\linebreak $(M^n,\bar g)$ with the induced
metric. Let $\xi$ be a vector field on $M^n$ along $F^l$. Then the
metric on $\xi(F^l)$, induced by the Sasaki metric of $TM^n$, is
defined by the following scalar product $$ {g_s}(\tilde X,\tilde
Y)=g\,(X,Y)+{\bar g}\,(\bar \nabla_X\, \xi ,\bar \nabla_Y\, \xi),
$$ for all vector fields $\tilde X=X^h+(\bar \nabla _X \, \xi)^v$
and $\tilde Y=Y^h+(\bar \nabla _Y \, \xi)^v$ on $\xi(F^l)$, where
$X,Y$ are vector fields on $F^l$.
\end{corollary}

\subsection{ Normal bundle of $\xi (F^l)$.}

To describe the normal bundle of $\xi(F^l)$, we need one auxiliary
notion. Let $\xi$ be a given vector field on a submanifold $F^l
\subset~M^n$. Then $\bar \nabla $ enables us to define a
point-wise linear mapping $\bar \nabla \xi: T_xF^l \to T_xM^n$, $X
\to \bar \nabla _X \xi$, for all $x \in M^n$. Its dual mapping,
with respect to the corresponding scalar products induced by $g$
and $\bar g$, gives rise to the linear mapping $(\bar \nabla
\xi)^*: T_xM^n \to T_xF^l$ defined by the formula

\begin{equation}\label{Eqn4}
g\,((\bar \nabla \xi)^*W,X)={\bar g}\,(\bar \nabla_X \xi,W) \mbox{
for all $W \in T_xM^n$ and $X \in T_xF^l$}.
\end{equation}
We call the mapping $(\bar \nabla \xi)^*: T_xM^n \to T_xF^l$ the
{\it conjugate derivative mapping}, or simply {\it conjugate
derivative}. Remark, that if $W$ is a vector field on $M^n$, then
the application of $(\bar\nabla\xi)^*$ gives rise to a vector
field $(\bar\nabla\xi)^*W$ on $F^l$  by
$
[(\bar \nabla\xi)^*W]_x=(\bar\nabla\xi)^*W_x\in T_xF^l \mbox{ for
all $x\in F^l$}.
$

Now we can prove

\begin{proposition}\label{Pr4}
Let $\eta$ and Z be normal and tangent vector fields on $F^l$
respectively. Then the lifts $$ \eta^h, \ \eta^v-((\bar \nabla
\xi)^* \eta)^h,\ Z^v-((\bar \nabla \xi)^* Z)^h $$ to the points of
$\xi(F^l)$ span the normal bundle of $\xi (F^l)$ in $TM^n$.
\end{proposition}
\begin{proof}
Let $\tilde X=X^h+(\bar \nabla_X \xi)^v$ be a vector field on $\xi
(F^l).$ Let $\eta $ and $Z$ be vector fields  on $F^l$ which are
normal and tangent to $F^l$ respectively. Taking into account
(\ref{Eqn1}), (\ref{Pr}) and (\ref{Eqn4}), we have $$
\begin{array}{l}
{g_s}(\tilde X, \eta^h)={\bar g}\,(X, \eta)=0 \\[2ex]
\begin{array}{rl}
{g_s}(\tilde X, \eta^v - [(\bar \nabla \xi)^* \eta]^h)= &-{\bar
g}\,(X,(\bar \nabla \xi)^*\eta) + {\bar g}\,(\bar \nabla_X \xi,
\eta) \\[1ex]
=&-{\bar g}\,(\bar \nabla_X \xi,\eta) + {\bar
g}\,(\bar \nabla_X \xi, \eta)=0
\end{array}\\[3ex]
\begin{array}{rl}
{g_s}(\tilde X, Z^v-[(\bar \nabla \xi)^*Z]^h)=&-{\bar g}\,(X,(\bar
\nabla \xi)^* Z) +{\bar g}\,(\bar \nabla_X \xi, Z) \\[1ex]
=&-{\bar g}\,(\bar \nabla_X \xi,Z)+{\bar g}\,(\bar \nabla_X
\xi,Z)=0
\end{array}
\end{array}
$$

Let $\eta _1, \dots , \eta_p$  ($p=1,\dots, n-l$) be a normal
frame of  $F^l$ while $f_1, \dots ,f_l$ span $T_xF^l$ at each
point $x\in F^l$. Consider the vector fields $$
N_\alpha=\eta^h_\alpha , \ \ P_\alpha=\eta^v_\alpha -((\bar \nabla
\xi)^* \eta_\alpha)^h,  \ \ F_i=f_i^v-((\bar \nabla \xi)^* e_i)^h,
$$ where $\alpha=1,\dots, n-l; \ i=1,\dots,n$. Let us show that
these are linearly independent. Indeed, suppose that $$
\lambda^\alpha N_\alpha + \mu^\alpha P_\alpha + \nu^iF_i=
\{\lambda^\alpha \eta_\alpha - \mu^\alpha(\bar \nabla
\xi)^*\eta_\alpha - \nu^i(\bar \nabla \xi)^*e_i\}^h+\{\mu^\alpha
\eta_\alpha + \nu^i f_i\}^v=0. $$ Because of the fact that the
horizontal and vertical components are linearly independent, we
see that $\mu^\alpha \eta_\alpha+\nu^if_i=0$ which is possible iff
$\mu^\alpha=0, \nu^i=0.$ Then, from the horizontal part of the
decomposition above we see that $\lambda^\alpha=0.$ So, $N_\alpha,
~P_\alpha$ and $F_i$ are linearly independent, which completes the
proof.

\end{proof}

\begin{remark}
In the case when $\xi $ is a normal vector field, the images
$(\bar\nabla \xi )^* \eta $ and $(\bar\nabla \xi )^*Z$ have a
simple and natural meaning, namely $$
\begin{array}{l}
(\bar\nabla \xi )^* \eta = g^{ik}{\bar g}\,(\nabla _k^\perp \xi,
\eta)\frac{\partial}{\partial u^i},
 \ \
(\bar\nabla \xi )^* Z = -A_\xi Z,
\end{array}
$$ where $ \nabla ^\perp $ is the normal bundle connection of
$F^l$  and $A_\xi$ is the shape operator of $F^l$ with respect to
the normal vector field $\xi$. In fact, $(\bar\nabla \xi )^* \eta$
is the vector field on $F^l$ dual to the 1-form ${\bar g}\,(\nabla
_k^\perp \xi, \eta)\,du^k$.
\end{remark}

\section{Characterization of submanifolds of $TM^n$ transverse to fibers.}

It is clear that all totally geodesic vector fields along
submanifolds of $M^n$ generate submanifolds in $TM^n$ which are
transverse to the fibers of $TM^n$. We study in this section the
converse question. We start with the local case.

\begin{proposition}\label{Transv}
Let $N^l$ be an embedded submanifold in the tangent bundle of a
Riemannian manifold $M^n$, which is transverse to the fiber at a
point $z\in N^l$, then there is a submanifold $F^l$ of $M^n$
containing $x=\pi(z)$, a neighborhood $U$ of $x$ in $M^n$, a
neighborhood $V$ of $z$ in $TM^n$ and a vector field $\xi$ on
$M^n$ along $F^l \cap U$ such that $N^l \cap V=\xi(F^l \cap U)$.
\end{proposition}
\begin{proof}
 Since $T_z N^l$ is transverse to the vertical subspace $V_z TM^n$ of
$TTM^n$ at $z$, $\pi_* \upharpoonright T_z N^l :T_z N^l \to T_x
M^n$ is injective, and so there is an open neighborhood $W$ of $z$
in $TM^n$ such that $\pi_* \upharpoonright T_{z'} N^l:T_{z'} N^l
\to T_{\pi(z')} M^n$ is injective for all $z'\in W \cap N^l$.
Hence $\pi \upharpoonright {W\cap N^l}:W\cap N^l \to M^n$ is an
immersion, and thus there exist a cubic centered coordinate system
$(U,\varphi)$ about $x=\pi(z)$ and a neighborhood $V$ of $z$ in
$W$ such that $\pi \upharpoonright {V\cap N^l}$ is 1:1 and $\pi(V
\cap N^l)$ is a part of a slice $F^l$ of $(U,\varphi)$ (\cite{Wr},
p. 28). The slice $F^l$ is a submanifold of $M^n$ and we have $\pi
\upharpoonright {V\cap N^l}: V\cap N^l \to U\cap F^l$ is an
imbedding onto, and so there is a $C^\infty$-mapping $\xi:F^l \cap
U \to N^l \cap V$ such that $\pi \circ \xi=Id_{F^l \cap U}$. In
other words, $\xi$ is a vector field on $M^n$ along $F^l \cap U$
such that $N^l \cap V=\xi(F^l \cap U)$.
\end{proof}
The global version of the last result requires further conditions.

\begin{theorem} \label{Trans-max}
Let $N^n$be a connected compact $n$-dimensional submanifold of the
tangent bundle of a connected simply connected Riemannian manifold
$M^n$, which is everywhere transverse to the fibers of $TM^n$.
Then $M^n$ is also compact, and there is a vector field $\xi$ on
$M^n$ such that $\xi(M^n)=N^n$.
\end{theorem}
\begin{proof}
The fact that $N^n$ is everywhere transverse to the fibers of
$TM^n$ implies that $\pi\upharpoonright {N^n}:N^n \to M^n$ is an
immersion. Since $M^n$ and $N^n$ are connected of the same
dimension and $N^n$ is compact, then $M^n$ is compact and
$\pi\upharpoonright {N^n}$ is a covering projection (cf.
\cite{KoNz}, Vol.~1, p.178). Now, $M^n$ is simply connected and so
$\pi \upharpoonright {N^n}$ is a diffeomorphism. Let $\xi:M^n \to
N^n$ be the inverse of $\pi \upharpoonright {N^n}$. Then $\xi$ is
a vector field on $M^n$ and $\xi(M^n)=N^n$.
\end{proof}

In a similar way, we can show the following:

\begin{theorem}
Let $N^l$ be a connected compact submanifold of the tangent bundle
of a connected simply connected manifold $M^n$, which is
transverse to the fibers it meets and projects onto a simply
connected submanifold $F^l$ of $M^n$. Then $F^l$ is compact and
there is a vector field $\xi$ on $M^n$ along $F^l$ such that
$\xi(F^l)=N^l$.
\end{theorem}

In the particular case of horizontal totally geodesic submanifolds
of $TM^n$, i.e. whose tangent space at any point is horizontal, we
can state the following:
\begin{theorem}\label{Hor-Sub}
Let $N^l$ be a connected complete totally geodesic horizontal
submanifold of the tangent bundle of a connected Riemannian
manifold $M^n$ which projects into a simply connected Riemannian
submanifold $F^l$ of $M^n$. Then $F^l$ is also complete and
totally geodesic in $M^n$ and there is a parallel vector field
$\xi$ on $M^n$ along $F^l$ such that $\xi(F^l)=N^l$.
\end{theorem}
\begin{proof}
By hypothesis, for all $z\in N^l$, $T_z N^l$ is a horizontal
subspace of $T_z TM^n$ with respect to the Levi-Civita connection
of $\bar g$. Hence $\pi \upharpoonright {N^l}:N^l \to F^l$ is an
isometric submersion of $N^l$ into $F^l$, with $N^l$ and $F^l$
connected and of the same dimension. Since $N^l$ is complete, also
$F^l$ is complete and $N^l$ is a covering space of $F^l$ (cf.
\cite{KoNz}, Vol.1, p.176). The fact that $F^l$ is simply
connected implies that $\pi \upharpoonright {N^l}:N^l \to F^l$ is
an isometry, and there is an isometry $\xi:F^l \to N^l$ such that
$\pi \upharpoonright {N^l}\circ \xi=Id_{F^l}$, i.e. $\xi$ is a
vector field on $M^n$ along $F^l$.

Now, $F^l$ is totally geodesic. Indeed, let $X$ and $Y$ be vector
fields on $F^l$, and denote by the same letters some of their
extensions to $M^n$. If we denote by $X^h$ and $Y^h$ their
horizontal lifts to $TM^n$, then $X^h \upharpoonright {N^l}$ and
$Y^h \upharpoonright {N^l}$ are vector fields on $TM^n$ along
$N^l$. For all $z\in N^l$, $T_z N^l$ being horizontal,  $\pi_*
\upharpoonright {T_z N^l}:T_z N^l \to T_x M^n$ is bijective. Since
$\pi_*(X^h(z))=X(\pi (z))$ and $\pi_*(Y^h(z))=Y(\pi (z))$, we have
that $X^h(z)$ and $Y^h(z)$ are tangent to $N^l$. Thus $(\tilde
\nabla_{X^h} Y^h)\upharpoonright {N^l}$ is tangent to $N^l$ and
hence horizontal. Consequently $(\tilde \nabla_{X^h} Y^h)
\upharpoonright {N^l}=(\bar \nabla_{X} Y)^h \upharpoonright {N^l}$
and is tangent to $N^l$. Hence $\bar \nabla_X Y=\pi_*\circ(\bar
\nabla_X Y)^h$ is tangent to $F^l$ and so $F^l$ is totally
geodesic. It remains to prove that $\xi$ is parallel along $F^l$.
In fact, for all $x\in F^l$ and $X\in T_x F^l$, the vector $X^h +
(\bar \nabla_X \xi)^v$ is tangent to $\xi(F^l)=N^l$ at $\xi(x)$
and is mapped onto $X$. Since $T_{\xi(x)}N^l$ is a horizontal
space, $\bar \nabla_X \xi=0$. Therefore, $\xi$ is parallel along
$F^l$.
\end{proof}

\begin{corollary}
Let $N^n$ be a connected complete totally geodesic horizontal
$n$-dimensional submanifold of the tangent bundle of a connected
simply connected Riemannian manifold $M^n$. Then $M^n$ is also
complete and there is a parallel vector field $\xi$ on $M^n$ such
that $\xi(M^n)=N^n$.
\end{corollary}

\section{ The conditions on $\xi (F^l)$ to be totally geodesic.}

Evidently, geometrical properties of the submanifold $\xi(F^l)$
depend on the submanifold $F^l$ and the vector field $\xi$. If one
does not pose any restrictions on them, the geometry of $\xi(F^l)$
becomes rather intricate. Nevertheless, it is possible to
formulate the conditions on $\xi(F^l)$ to be totally geodesic in
more or less geometrical terms.

To do this, we introduce the notion of a $\xi$-connection on the
Riemannian manifold $M^n$.
\begin{definition}
Let $M^n$ be  a Riemannian manifold with Riemannian connection
$\bar\nabla$ and curvature tensor $\bar R$. Let $\xi$ be a fixed
smooth vector field on $M^n$. Denote by $\mathfrak{X}(M^n)$ the
set of all smooth vector fields  on $M^n$. The mapping
$\stackrel{*}{\nabla}:\mathfrak{X}(M^n)\times \mathfrak{X}(M^n)\to
\mathfrak{X}(M^n)$ defined by

\begin{equation}\label{Conn}
\nablast_{\bar X}{\bar Y}=\bar\nabla_{\bar X}\bar Y+\frac12\Big[ \bar R(\xi,\bar\nabla_{\bar
X}\xi)\bar Y+ \bar R(\xi,\bar\nabla_{\bar Y}\xi)\bar X\Big]
\end{equation}
is a torsion-free affine connection on $M^n$. It is called the $\xi$-connection.
\end{definition}

Remark that if $\xi$ is a parallel vector field  or the manifold $M^n$ is flat, then
the $\xi$-connection is the same as the Levi-Civita connection of $M^n$.

It is easy to check that (\ref{Conn}) indeed defines a torsion-free affine connection.
Now we can state the main technical tool for the further considerations.

\begin{proposition}\label{Pr5}
Let $F^l$ be a submanifold in a Riemannian manifold $M^n.$ Let $\xi$ be a vector field
on  $M^n$ along $F^l$. Then $\xi(F^l)$ is totally geodesic in $TM^n$ if and only if
\begin{itemize}
\item[(a)] $F^l$ is totally geodesic with respect to
the $\xi$-connection (\ref{Conn});
\item[(b)] for any vector fields $X,Y$ on $F^l$
$$
\bar \nabla_X \bar \nabla_Y \xi =
\bar \nabla_{\nablast_X Y} \xi +\frac{1}{2} \bar R(X,Y)\xi.
$$
\end{itemize}
\end{proposition}
\begin{proof}
By definition, the submanifold $\xi(F^l)$ is totally geodesic in $TM^n$  if and only
if $g_s\,(~\tilde \nabla_{\tilde X}\tilde Y,\tilde N)~=~0$ for any vector fields
$\tilde X,\tilde Y$ tangent to $\xi(F^l)$ along $\xi(F^l)$ and $\tilde N$ normal to
$\xi (F^l)$. To calculate $\tilde \nabla_{\tilde X}\tilde Y$, we use  the  Kowalski
formulas \cite{Kow}.

{\it For any vector fields $\bar X,\bar Y$ on $M^n$, the covariant derivatives of various
combinations of lifts to the point $(x,\xi) \in TM^n$ can be found as follows}

\begin{equation}\label{Kow}
\begin{array}{ll}
\tilde \nabla_{\bar X^h}\bar Y^h = (\bar \nabla_{\bar X} \bar Y)^h-
\frac{1}{2}(\bar R (\bar X,\bar Y) \xi)^v, \
&\tilde \nabla_{\bar X^v}\bar Y^h = \frac{1}{2} (\bar R (\xi ,\bar X) \bar Y)^h,\\[2ex]
\tilde \nabla_{\bar X^h}\bar Y^v = (\bar \nabla_{\bar X} \bar Y)^v+
\frac{1}{2}(\bar R (\xi ,\bar Y) \bar X)^h, \
& \tilde \nabla_{\bar X^v}\bar Y^v = 0.
\end{array}
\end{equation}
{\it where $\bar \nabla$ and $\bar R$ are the Levi-Civita connection and the curvature
tensor of $M^n$ respectively}.

Let $\tilde X=X^h+(\bar \nabla_X \xi)^v$ and $\tilde Y=(Y)^h+(\bar \nabla_Y \xi)^v$ be
vector fields tangent to $\xi(F^l).$ Then, applying (\ref{Kow}), we easily
find
$$
\tilde \nabla_{\tilde X} \tilde Y = (\bar \nabla_X Y+\frac{1}{2} \bar R(\xi ,\bar
\nabla _X \xi)Y+\frac{1}{2}\bar R(\xi,\bar\nabla_Y\xi)X)^h+ (\bar \nabla_X \bar
\nabla_Y\, \xi - \frac{1}{2}\bar R(X,Y)\xi)^v
$$
or
$$
\tilde \nabla_{\tilde X} \tilde Y = (\nablast_X Y)^h+ (\bar \nabla_X \bar
\nabla_Y\, \xi - \frac{1}{2}\bar R(X,Y)\xi)^v.
$$

Using Proposition \ref{Pr4}, we see that the totally geodesic property of $\xi(F^l)$
is equivalent to
\begin{equation}\label{Cond}
\left\{
\begin{array}{rl}
{\bar g}\,(\nablast_X Y ,\eta)&=0,\\[2ex]
{\bar g}\,(\nablast_X Y,(\nabla \xi)^* \eta)&={\bar g}\,(\bar \nabla_X \bar \nabla_Y
\xi-\frac{1}{2}\bar R(X,Y)\xi ,\eta),\\[2ex]
{\bar g}\,(\nablast_X Y,(\nabla \xi)^*Z)&={\bar g}\,(\bar \nabla_X \bar \nabla_Y
\xi-\frac{1}{2}\bar R(X,Y)\xi ,Z),
\end{array}
\right.
\end{equation}
for any vector fields $X,Y,Z$ tangent to $F^l$ and any vector field $\eta$ orthogonal
to $F^l$.

From $(\ref{Cond})_1$ we see that $F^l$ must be autoparallel with respect to
$\nablast$ and hence totally geodesic \cite{KoNz}. Thus, $\nablast_XY$ is tangent to
$F^l$ and it is possible to apply (\ref{Eqn4}). Therefore, we can rewrite the
equations $(\ref{Cond})_2$ and $(\ref{Cond})_3$ as
$$
\left\{
\begin{array}{l}
{\bar g}\,(\bar \nabla_{\nablast_X Y} \xi - \bar \nabla_X \bar \nabla_Y \xi+
\frac{1}{2} \bar R(X,Y) \xi,\eta) =0, \\[1ex]
{\bar g}\,(\bar \nabla_{\nablast_X Y} \xi -\bar \nabla_X \bar \nabla_Y \xi +
\frac{1}{2}\bar R(X,Y) \xi,Z) =0
\end{array}
\right.
$$
for any vector fields $\eta $ normal and $Z$ tangent to $F^l$ along $F^l$. Thus, we conclude
$$
\bar \nabla_X \bar \nabla_Y \xi =\bar \nabla_{\nablast_X Y} \xi
+ \frac{1}{2} \bar R(X,Y) \xi,
$$
which completes the proof.
\end{proof}

For the cases when $l=1$ and $l=n$, we get the known conditions for the totally
geodesic property of $\xi(F^l)$.
\begin{corollary}\label{l=1}
If $l=1$ and $\xi(F^l)$  is a curve $\Gamma$ in $TM^n$ then
this curve is a geodesic if and only if
$$
\left\{
\begin{array}{l}
 x''+\bar R(\xi,\xi')x'=0, \\[1ex]
 \xi''=0,
\end{array}
\right.
$$
where $(')$ means the covariant derivative with respect to the natural parameter of
$\Gamma$ and $x(\sigma)=(\pi\circ\Gamma)(\sigma)$ \emph{(cf. \cite{Sk})};
\end{corollary}

\begin{proof}
Indeed, in this case $\tilde X=\tilde Y=\Gamma'=(x')^h+(\xi')^v$, $\bar X=\bar Y=x'$
and $\nablast_{\bar X}{\bar Y}=x''+\bar R(\xi,\xi')x'$. Thus, $x(\sigma)$ is geodesic
with respect to the $\xi$-connection iff $x''+\bar R(\xi,\xi')x'=0$ and the rest of
the proof is evident.
\end{proof}

\begin{corollary}\label{l=n}
 If $l=n$ and $F^l=M^n$,  then $\xi(M^n)$ is totally geodesic in $TM^n$ if and only if
 for any vector fields $\bar X,\bar Y$ on $M^n$ \emph{(cf. \cite{Wk})}
$$
\bar \nabla_{\bar X} \bar \nabla_{\bar Y} \xi =
\bar \nabla_{\nablast_{\bar X}\bar Y} \xi +\frac{1}{2} \bar R(\bar X,\bar Y)\xi.
$$
\end{corollary}
\begin{proof}
In this case, only $(b)$ of Proposition \ref{Pr5} should be checked, which completes
the proof.
\end{proof}

The result of Corollary \ref{l=n} can be expressed in more
geometrical terms. To do this, introduce a symmetric bilinear
mapping $h_\xi:  \mathfrak{X}(M^n)\times \mathfrak{X}(M^n)\to
\mathfrak{X}(M^n)$ by
\begin{equation}\label{h}
h_\xi(\bar X,\bar Y)=\frac12 \Big[\bar R(\xi,\nabla_{\bar X}\xi)\bar Y+
     \bar R(\xi,\nabla_{\bar Y}\xi)\bar X\Big],
\end{equation}
for all $\bar X$, $\bar Y \in \mathfrak{X}(M^n)$. Then the
definition of the $\xi$-connection takes as similar form as the
Gauss decomposition
\begin{equation}\label{Conn1}
\nablast_{\bar X}{\bar Y}=\bar\nabla_{\bar X}{\bar Y}+h_\xi(\bar X,\bar Y).
\end{equation}
Define a \textit{``shape operator"} $A_\xi$ for the field $\xi$ by
\begin{equation}\label{Shp}
A_\xi\bar Y=-\bar \nabla_{\bar Y}\xi,\;\textup{for all}\; \bar Y
\in \mathfrak{X}(M^n).
\end{equation}
Then the covariant derivative of the $(1,1)$-tensor field $A_\xi$ is given by
$$
(\bar\nabla_{\bar X}A_\xi)\bar Y=-\bar\nabla_{\bar X}\bar \nabla_{\bar
Y}\xi+\bar\nabla_{\bar\nabla_{\bar X}\bar Y}\xi.
$$
Hence we see that the Codazzi-type equation
$
\bar R(\bar X,\bar Y)\xi=(\bar\nabla_{\bar Y}A_\xi)\bar X-(\bar\nabla_{\bar X}A_\xi)\bar Y
$
holds. In these notations $$ \bar\nabla_{\nablast_{\bar X}\bar Y}
\xi +\frac{1}{2} \bar R(\bar X,\bar Y)\xi-\bar \nabla_{\bar X}
\bar \nabla_{\bar Y} \xi= \bar\nabla_{h_\xi(\bar X,\bar Y)}\xi+
\frac{1}{2}\Big[(\bar\nabla_{\bar X}A_\xi)\bar Y+(\bar\nabla_{\bar
Y}A_\xi)\bar X\Big]. $$ If we introduce a symmetric bilinear
mapping $ \Omega _\xi:  \mathfrak{X}(M^n)\times
\mathfrak{X}(M^n)\to \mathfrak{X}(M^n)$ defined by $$ \Omega _\xi
(\bar X,\bar Y)=\bar\nabla_{h_\xi(\bar X,\bar Y)}\xi+
\frac{1}{2}\Big[(\bar\nabla_{\bar X}A_\xi)\bar Y+(\bar\nabla_{\bar
Y}A_\xi)\bar X\Big], $$ then Corollary \ref{l=n} can be
reformulated as
\begin{corollary}
 If $\xi$ is a smooth  vector field on a Riemannian manifold $M^n$ then
 $\xi(M^n)$ is totally geodesic in $TM^n$ if and only if for any vector fields
 $\bar X,\bar Y$ on $M^n$
\begin{equation}\label{Omega}
\Omega _\xi (\bar X,\bar Y)=\bar\nabla_{h_\xi(\bar X,\bar Y)}\xi+
\frac{1}{2}\Big[(\bar\nabla_{\bar X}A_\xi)\bar Y+(\bar\nabla_{\bar Y}A_\xi)\bar
X\Big]\equiv 0,
\end{equation}
where $h_\xi$ and $A_\xi$ are defined by (\ref{h}) and (\ref{Shp}) respectively.
\end{corollary}

\begin{remark} The statement of Proposition \ref{Pr5} can also be reformulated in
these terms, namely, {\it let $F^l$ be a  submanifold in a Riemannian manifold $M^n$
and $\xi $ be a vector field on $M^n$ along $F^l$. Then $\xi (F^l)$ is totally
geodesic in $TM^n$ if and only if $F^l$ is totally geodesic with respect to the
$\xi$-connection (\ref{Conn}) and $\Omega_\xi$ vanishes on the tangent bundle of
$F^l$}
\end{remark}

Now, combining Theorem \ref{Trans-max} with Proposition \ref{Pr5},
we obtain

\begin{corollary}
On a connected simply connected compact $n-$dimensional Riemannian
manifold, vector fields satisfying $(b)$ of Proposition \ref{Pr5}
generate the only connected compact totally geodesic
$n$-dimensional submanifolds of the tangent bundle which are
transverse to fibers.
\end{corollary}

As has been shown in \cite{Ym}, for the case of the unit tangent bundle, the Hopf
vector fields on odd dimensional spheres generate totally geodesic submanifolds in
$T_1S^n$. For the tangent bundle the situation is different.

\begin{theorem}
 A non-zero Killing vector field  on a space of non-zero constant curvature $(M^n,c)$ never
generates a  totally geodesic submanifold in $TM^n$. Moreover, a manifold with
positive sectional curvature does not admit a non-zero Killing vector field with
totally geodesic property.
\end{theorem}
\begin{proof}
Let $\xi$ be a Killing vector field on a space $M^n$ of constant curvature $c$. Then
$A_\xi $ is a skew-symmetric linear operator, i.e.
\begin{equation}\label{Killing}
 \bar g(A_\xi \bar X,\bar Y)+\bar
g(\bar X,A_\xi \bar Y)=0,
\end{equation}
and moreover,
\begin{equation}\label{KillProp}
(\bar\nabla_{\bar X}A_\xi)\bar Y=\bar R(\xi,\bar X)\bar Y
\end{equation}
for all vector fields $\bar X,\bar Y$ on $M^n$ (cf. \cite{KoNz}). Since $M^n$ is of
non-zero constant curvature, the equation (\ref{Omega}) can be simplified in the
following way.
$$
\begin{array}{rl}
(\bar\nabla_{\bar X}A_\xi)\bar Y+(\bar\nabla_{\bar Y}A_\xi)\bar X=&\bar R(\xi,\bar X)\bar Y+\bar R(\xi,\bar Y)\bar
X=\\[1ex]
&c\,\Big[2\bar g(\bar X,\bar Y)\,\xi-\bar g(\xi,\bar X)\bar Y-\bar g(\xi,\bar Y)\bar
X\Big]
\end{array}
$$
$$
\begin{array}{l}
\bar R(\xi,\bar\nabla_{\bar X}\xi)\bar Y+\bar R(\xi,\bar\nabla_{\bar Y}\xi)\bar X=
c\,\Big[\bar g(\bar\nabla_{\bar X}\xi,\bar Y)+\bar g(\bar X,\bar\nabla_{\bar Y}\xi)\bar X)\Big]\xi-\\[1ex]
c\,\Big[(\bar g(\xi,\bar X)\bar\nabla_{\bar Y}\xi+\bar g(\xi,\bar Y)\bar\nabla_{\bar
X}\xi)\Big]= c\,\Big[\bar g(\xi,\bar X)A_\xi\bar Y+\bar g(\xi,\bar Y)A_\xi\bar X\Big].
\end{array}
$$
So, $\xi$ is totally geodesic if
$$
 \bar g(\xi,\bar X)\bar Y+\bar g(\xi,\bar Y)\bar X
 -\bar\nabla_{\bar g(\xi,\bar X)A_\xi\bar Y+\bar g(\xi,\bar Y)A_\xi\bar X}\xi=
 2\bar g(\bar X,\bar Y)\,\xi,
$$
or
$$
 \bar g(\xi,\bar X)\Big[\bar Y+A_\xi(A_\xi\bar Y)\Big]+\bar g(\xi,\bar Y)\Big[\bar X
 +A_\xi(A_\xi\bar X)\Big]=2\bar g(\bar X,\bar Y)\,\xi,
$$
for all vector fields $\bar X,\bar Y$ on $M^n$. Choosing $\bar X,\bar Y$ such that
$\bar X_x\ne 0$ and $\bar X_x=\bar Y_x \perp \xi_x,$ we get $2|\bar X_x|^2\xi_x=0$.
Therefore, $\xi=0$ for all $x\in M^n.$

Let $\xi$  be a non-zero Killing vector field on a manifold with  \textit{positive}
(non-constant) sectional curvature.  From (\ref{Killing}) it follows that
$A_\xi\xi\perp\xi$. If $A_\xi\xi=0$, then, after setting $Y=\xi$ in (\ref{Killing}),
we conclude that $\xi$ has a constant length and therefore can be totally geodesic if
it is a parallel vector field \cite{Wk}. In this case, $M^n=M^{n-1}\times E^1$ and we
come to a contradiction. Suppose that $A_\xi\xi\ne0$. Then $\xi\wedge A_\xi\xi$ is a
non-zero bivector field. Setting $\bar Y=\bar X$ in (\ref{Omega}) and using
(\ref{KillProp}), we have
$$
A_\xi\Big[\bar R(\xi,A_\xi \bar X)\bar X\Big]+\bar R(\xi, \bar X)\bar X=0.
$$
Taking a scalar product in both sides with $\xi$ and applying (\ref{Killing}), we get
$$
-\bar g(\bar R(\xi,A_\xi \bar X)\bar X,A_\xi\xi)+K_{\xi\wedge \bar X}|\xi\wedge \bar
X|^2=0.
$$
Finally, setting $\bar X=A_\xi\xi$, we have $K_{\xi\wedge \bar X}=0$ and come to a
contradiction.

\end{proof}

The next Theorem is analogous to the one proved by Walczak P. \cite{Wk}, but does not
have similar rigid consequences for the structure of $M^n$.
\begin{theorem}\label{Th1}
Let $\xi$ be a vector field of constant length along a submanifold $F^l \subset M^n$.
Then $\xi(F^l)$ is a totally geodesic submanifold in $TM^n$ if and only if $F^l$ is
totally geodesic in $M^n$ and $\xi$ is a parallel vector field on $M^n$ along $F^l$.
\end{theorem}
\begin{proof}
 The condition $|\,\xi\, |=const$ implies ${\bar g}\,(\bar \nabla_X \xi, \xi)=0$ for
any vector field $X$ tangent to $F^l$ . As $\xi(F^l)$ is supposed to be totally
geodesic, it follows from the second condition of Proposition \ref{Pr5} that ${\bar
g}\,(\bar \nabla_X \bar \nabla_Y \xi, \xi)=0$. Hence ${\bar g}\,(\bar \nabla_X \xi,
\bar \nabla_Y \xi)=0$ for any $X,Y \in T_xF^l$, $x \in F^l$. Supposing $X=Y$, we see
that $\bar \nabla_X \xi =0$, i.e. $\xi$ is parallel along $F^l$ in the ambient space
and the second condition of Proposition \ref{Pr5} is fulfilled. Moreover, the
condition $\bar \nabla_X \xi =0$ means that the $\xi$-connection (\ref{Conn})
coincides with the Levi-Civita connection of $M^n$, so that by Proposition \ref{Pr5}
$F^l$ is totally geodesic in $M^n$.

On the other hand, if $F^l$ is totally geodesic in $M^n$ and $\bar \nabla_X \xi=0$ for
any tangent vector field $X$ on $F^l$, then both conditions from Proposition \ref{Pr5}
are satisfied evidently.

\end{proof}

Giving more restrictions on the vector field, we can a more geometrical result.
\begin{theorem}\label{Th2}
Let $\xi $ be  a normal vector field  on a submanifold $~F^l \subset~M^n,$ which is
parallel in the normal bundle. Then $\xi (F^l)$ is totally geodesic in $TM^n$ if and
only if $F^l$ is totally geodesic in $M^n.$
\end{theorem}
\begin{proof}
If $\xi$ is a normal vector field to $F^l$ and parallel in the normal bundle, then
$\bar \nabla_X \xi=-A_{\xi} X$ for each vector field $X$ on $F^l$, where $A_{\xi}$ is
the shape operator  of $F^l$ with respect to $\xi,$ and hence ${\bar g}\,(\bar
\nabla_X \xi, \xi)=0.$ This means that $|\xi |$=const along $F^l$.

Let $\xi(F^l)$ be totally geodesic in $TM^n$. Then from (b) of Proposition \ref{Pr5}
we see that $\bar g \,(\bar\nabla_X\bar\nabla_Y\xi,\xi)=0$, which implies $|\bar
\nabla_X\xi|=0$ for each $X$ tangent to $F^l$.  In this case, along $F^l$ the
$\xi$-connection (\ref{Conn}) coincides with the Levi-Civita connection of $M^n$ and
(a) of Proposition \ref{Pr5} implies the totally geodesic property of $F^l$.

Conversely, if $\xi$ is a normal vector field which is parallel in the normal bundle
of $F^l$ and $F^l$ is totally geodesic, then $\bar \nabla_X\xi=0$ for any vector field
$X$ tangent to $F^l$. Evidently, both conditions of Proposition \ref{Pr5} are
fulfilled.
\end{proof}

 The application of Theorem \ref{Th2} to the specific case of  a foliated Riemannian
manifold allows to clarify the geometrical structure of $\xi(M^n)$. The manifold $M^n$
is said to be \textit{$\nu$-foliated} if it admits a family $\mathcal{F}$ of connected
$\nu$-dimensional submanifolds $\{\mathcal{F}_\alpha; \alpha\in A\}$ called
\textit{leaves} such that (i) $M^n=\bigcup\limits_{\alpha\in A}\mathcal{F}_\alpha$;
(ii) $\mathcal{F}_\alpha\cap\mathcal{F}_\beta=\emptyset$ for $\alpha\ne\beta$; (iii)
there exists a coordinate covering  $\mathcal{U}$ of $M^n$ such that in each local
chart $U\in\mathcal{U}$ the leaves can be expressed locally as level submanifolds,
i.e. $u^{\nu+1}=c_{\nu+1},\dots,u^n=c_n$.

 The family $\mathcal{F}$ is called a \textit{$\nu$-foliation} and
\textit{hyperfoliation} for $\nu=n-1$. The hyperfoliation is said to be
\textit{transversally orientable }if $M^n$ admits a vector field $\xi$ transversal to
the leaves. Moreover, with respect to the Riemannian metric on $M^n$, this vector
field can be chosen as a field of unit normals for each leaf.

A submanifold $F^{k+\nu}\subset M^n$ is called \textit{$\nu$-ruled} if $F^{k+\nu}$
admits a $\nu$-foliation $\big\{\mathcal{F}_\alpha;\, \alpha\in A\big\}$ such that
each leaf $\mathcal{F}_\alpha$ is totally  geodesic in $M^n$. The leaves
$\mathcal{F}_\alpha$ are called \textit{elements} or \textit{generators} \cite{Rov}.

\begin{corollary}\label{Foli}
Let $M^n$ be a Riemannian manifold admitting a totally geodesic transversally
orientable hyperfoliation $\mathcal{F}$. Let $\xi$ be a field of normals of the
foliation having constant length.  Then $\xi(M^n)$ is an $(n-1)$-ruled submanifold in
$TM^n$ with the elements $\xi(\mathcal{F}_\alpha)$.
\end{corollary}
\begin{proof}
Indeed, let $\mathcal{F}_\alpha$ be a leaf of the hyperfoliation and $\xi$ be a vector
field of constant length on $M^n$ which is a field of normals along each leaf.
Applying Theorem \ref{Th2}, we get that $\xi(\mathcal{F}_\alpha)$ is totally geodesic
in $TM^n$ for each $\alpha$. Since $\xi:M^n\to \xi(M^n)$ is a homeomorphism,
$\xi(\mathcal{F}_\alpha)\cap\xi(\mathcal{F}_\beta)=\emptyset$ for $\alpha\ne\beta$ and
$\xi(M^n)=\bigcup\limits_{\alpha\in A}\xi(\mathcal{F}_\alpha)$. Finally, if
$\mathcal{F}_\alpha$ is given by $u^{n}=c_n$ within a local chart $U$ then from
(\ref{Sub}) we see that $\xi(\mathcal{F}_\alpha)$ is given by the same equalities
within the local chart $\xi(U)$. So, $\xi(\mathcal{F})=\big\{\xi(\mathcal{F}_\alpha);
\alpha\in A\big\}$ form a hyperfoliation on $\xi(M^n)$ with totally geodesic leaves in
$TM^n$.
\end{proof}

\section{The case of a base space of constant curvature.}

If the ambient space is of constant curvature $c\ne0$ and $\xi$ is a normal vector
field on a submanifold $F^l \subset M^n$, then the necessary and sufficient condition
on $\xi$ to generate a totally geodesic submanifold in $TM^n$ takes a rather simple
form.
\begin{theorem}\label{Th3}
Let $F^l$ be a submanifold of a space $M^n(c)$ of constant curvature $c\ne 0$. Let
$\xi$ be a normal vector field on $F^l.$ Then $\xi (F^l)$ is totally geodesic in
$TM^n$ if and only if $F^l$ is totally geodesic in $M^n(c)$ and $\xi$ is parallel in
the normal bundle.
\end{theorem}
\begin{proof}
The curvature tensor of $M^n (c)$ is of the form
\begin{equation}\label{R}
\bar R(\bar X,\bar Y)\bar Z=
c\ ( {\bar g}\,(\bar Y,\bar Z)\bar X-{\bar g}\,(\bar X,\bar Z)\bar Y\,).
\end{equation}
If $\xi $ is a normal vector field on $F^l$ then
$
\bar \nabla_X \xi=-A_{\xi} X+\nabla^{\perp}_X \xi.
$
As $A_{\xi} X$ is tangent and $\nabla^{\perp}_X \xi$
is normal to $F^l$, from (\ref{R}) we find
$$
\bar R(\xi, \bar \nabla_X \xi) Y= -c\,g\,(A_{\xi} X,Y)\,\xi
$$
for any vector fields $X,Y$ on $F^l.$ Thus, the conditions from Proposition \ref{Pr5}
mean that
\begin{equation}\label{Eqn5}
\left\{
\begin{array}{l}
\bar \nabla_X Y-c\,g\,(A_{\xi} X,Y) \xi
\mbox{ \ \ is tangent to \ } F^l, \\[2ex]
\bar \nabla_{\bar \nabla_X Y-c\,g\,(A_{\xi} X,Y) \xi} \xi = \bar \nabla_X
\bar \nabla_Y \xi.
\end{array}
\right.
\end{equation}

Multiplying $(\ref{Eqn5})_1$ by $\xi$ and by normal vector field $\eta$ orthogonal to
$\xi$, we have
$$
\left\{
\begin{array}{r}
g\,(A_{\xi} X,Y)(1-c \, |\xi |^2)=0,\\[1ex]
g\,(A_{\eta} X,Y)=0.
\end{array}
\right.
$$

If $\xi$ is of constant length $|\xi |^2=\frac{1}{c} \ \ (c>0)$
then by Theorem \ref{Th1}, $F^l$ is totally geodesic in $M^n,$ otherwise $F^l$
is totally geodesic immediately.

So, $F^l$ is totally geodesic and therefore $\bar \nabla_X \xi= \nabla^{\perp}_X \xi$,
$\bar\nabla_XY=\nabla_XY$. The condition $(\ref{Eqn5})_2$ now takes the form
\begin{equation}\label{Eqn6}
\nabla^{\perp}_{\nabla_X Y} \xi = \nabla^{\perp}_X
\nabla^{\perp}_Y \xi.
\end{equation}

Set $Y=\nabla_V Z$, where $V$ and $Z$ are arbitrary vector fields tangent to $F^l$.
Then from (\ref{Eqn6}), we get
$$
\nabla^{\perp}_{\nabla_X \nabla_V Z} \xi = \nabla^{\perp}_X \nabla^
{\perp}_{\nabla_V Z} \xi.
$$
Applying (\ref{Eqn6}) to $\nabla^{\perp}_{\nabla_V Z} \xi$ in the right-hand side of
the above equation, we see that $\nabla^{\perp}_{\nabla_V Z} \xi =
\nabla^{\perp}_V\nabla^{\perp}_Z \xi$ and therefore,

\begin{equation}\label{Eqn7}
\nabla^{\perp}_{\nabla_X \nabla_V Z} \xi = \nabla^{\perp}_X \nabla^{\perp}_V
\nabla^{\perp}_Z \xi.
\end{equation}
Interchanging the roles of $X$ and $V$, we get

\begin{equation}\label{Eqn8}
\nabla^{\perp}_{\nabla_V \nabla_X Z} \xi = \nabla^{\perp}_V \nabla^{\perp}_X
\nabla^{\perp}_Z \xi.
\end{equation}
Finally, applying again (\ref{Eqn6}) to the bracket $[X,V]$ and $Z$, we get

\begin{equation}\label{Eqn9}
\nabla^{\perp}_{\nabla_{[X,V]}Z} \xi = \nabla^{\perp}_{[X,V]}
\nabla^{\perp}_Z \xi.
\end{equation}

Combining (\ref{Eqn7}),(\ref{Eqn8}) and (\ref{Eqn9}), we obtain
$$
\nabla^{\perp}_{R(X,V)Z} \xi = R^\perp(X,V) \nabla^{\perp}_Z \xi
$$
where $R$ is the curvature tensor of $F^l$ and $R^\perp$ is the normal curvature
tensor. Since $F^l$ is totally geodesic and $M^n(c)$ is of constant curvature,
$R^\perp(X,Y)\eta \equiv 0$ for any normal vector field $\eta$ and, moreover,
$$
R(X,Y)Z = c\,(g\,(Y,Z)X - g\,(X,Z)Y).
$$
So, we have
$$
c\, \nabla^{\perp}_{g\,(Y,Z)X-g\,(X,Z)Y} \xi = 0.
$$
Setting $X$ orthogonal to $Y$ and $Y=Z$ we get
$
\nabla^{\perp}_X \xi = 0
$
for any vector field $X$ on $F^l$, which completes the necessary part of the proof.
The sufficient part is trivial.
\end{proof}

The application of Theorem \ref{Th3} to the case of a space of
constant curvature shows the difference between our considerations
and Walczak's \cite{Wk}. Let $S^n$ be the unit sphere and
$S^{n-1}$ be the unit totally geodesic great sphere in $S^n$.
Denote by $D^n$ an open equatorial zone around $S^{n-1}$ where the
unit geodesic vector field orthogonal to $S^{n-1}$ is regularly
defined. Then $D^n$ is a Riemannian manifold of constant positive
curvature and $S^{n-1}$ is a totally geodesic submanifold in
$D^n$.

\textit{Let $\xi$ be a unit (or of constant length) geodesic
vector field on $D^n\subset S^n$ which is normal to the totally
geodesic great sphere $S^{n-1}$. Then $\xi(D^n)$ is not totally
geodesic in $TD^n$ while the restriction of $\xi$ to $S^{n-1}$
generates the totally geodesic submanifold $\xi(S^{n-1})$ in
$TD^n$.}

Indeed, $\xi$ is of constant length and by Walczak's result,
$\xi(D^n)$  can be totally geodesic in $TD^n$ only if $\xi$ is a
parallel vector field on $D^n$ \cite{Wk}, which is impossible due
to positive curvature of $D^n$. On the other hand, $\xi$ is
parallel in the normal bundle of $S^{n-1}\subset D^n$ and we can
apply Theorem \ref{Th3} to see that $\xi(S^{n-1})$ is totally
geodesic in $TD^n$.

%%%%%%%%%%%%%%%%%%%%%%%%%%%%%%%%%%%%%%%%%%%%%%%
As concerns flat Riemannian manifolds, Walczak has shown that
every totally geodesic vector field on a flat Riemannian manifold
is harmonic (cf. \cite{Wk}) and that, consequently, on a compact
flat Riemannian manifold, a vector field is totally geodesic if
and only if it is parallel. We shall give a similar result for
vector fields along submanifolds.

\begin{theorem}\label{Flat}
Let $F^l$ be a compact oriented submanifold in a flat Riemann\-ian
manifold $M^n$. Let $\xi$ be a vector field on $F^l$. Then
$\xi(F^l)$ is totally geodesic in $TM^n$ if and only if $F^l$ is
totally geodesic in $M^n$ and $\xi$ is parallel along $F^l$.
\end{theorem}
\begin{proof}
Since $M^n$ is flat, the $\xi$-connection is the same as the
Levi-Civita connection on $M^n$. So, by Proposition \ref{Pr5},
$\xi(F^l)$ is totally geodesic if and only if $F^l$ is totally
geodesic and
\begin{equation}\label{A19}
\bar \nabla_{X}\bar \nabla_{Y}\xi=\bar \nabla_{\bar \nabla_{X}Y} \xi
\end{equation}
for all vector fields $X$ and $Y$ on $F^l$.

Suppose now that $\xi(F^l)$ is totally geodesic. Then $F^l$ is
totally geodesic and is thus flat. Hence locally we can choose
vector fields $X_1$, $X_2$,...,$X_l$ tangent to $F^l$ such that
$\bar \nabla_{X_i}X_j=\nabla_{X_i}X_j=0$, and $\bar
g(X_i,X_j)=g(X_i,X_j)=\delta_{ij}$, for all $i,j=1,...,l$. Putting
$X=Y=X_i$ in the identity (\ref{A19}), we obtain $\bar
\nabla_{X_i}\bar \nabla_{X_i}\xi=0$. Hence, $\sum_{i=1}^l \bar
g(\bar \nabla_{X_i}\bar \nabla_{X_i}\xi,\xi)=0$, i.e.

\begin{equation}\label{A20}
\sum_{i=1}^l X_i.\bar g( \bar \nabla_{X_i}\xi,\xi)=\sum_{i=1}^l |
\bar\nabla_{X_i}\xi |^2.
\end{equation}

If we consider the function $f$ defined by $f(x)=\frac {1}{2} \bar
{g}_x(\xi,\xi)$, for all $x\in F^l$, then we can define a global
vector field $X_f$ on $F^l$ by the local formula $X_f =g( \bar
\nabla_{X_i}\xi,\xi)X_i$. Formula (\ref{A20}) can thus be written
locally as div$X_f$=$\sum_{i=1}^l | \bar \nabla_{X_i}\xi |^2$.

Integrating both sides of the last equality and applying Green's
theorem, we obtain $\sum_{i=1}^l \int_{F^l}| \bar \nabla_{X_i}\xi
|^2 dv=0$, and hence $\bar\nabla_{X_i}\xi=0$, for all $i=1,...,l$.
Therefore $\xi$ is parallel along $F^l$.

The sufficient part of the theorem is trivial.

\end{proof}
\begin{remarks}

1.\ If in Theorem \ref{Flat} the field $\xi$ is a normal vector field along $F^l$,
then $\bar \nabla_{X}\xi$ is also normal for each vector field $X$ on $F^l$. Indeed,
for the $X_i$'s constructed in the proof of the theorem, we have $\bar g( \bar
\nabla_{X_i}\xi,X_j)=X_i.\bar g(\xi,X_j)=0$, and so $ \bar \nabla_{X_i}\xi$ is normal
to $F^l$. Hence the identity (\ref{A19}) can be written as

\begin{equation}\label{A21}
\bar \nabla_{X}^{\perp}\bar \nabla_{Y}^{\perp}\xi=\bar
\nabla_{\nabla_X Y}\xi.
\end{equation}

Also, $\xi$ is parallel if and only if it is parallel in the normal bundle. Hence
$\xi(F^l)$ is totally geodesic if and only if $F^l$ is totally geodesic and $\xi$ is
parallel in the normal bundle.

2.\ The condition of compactness is necessary. Indeed, if we consider $\mathbb{R}^n$
with its canonical coordinates $(x_1,x_2,...,x_n)$ and its canonical Euclidean metric,
and the hypersurface $\mathbb{R} ^{n-1}$ which is identified with the subspace given
by: $x_n=0$, then $\mathbb{R} ^{n-1}$ is an oriented totally geodesic submanifold of
$\mathbb{R} ^n$. We have $\bar \nabla_{\partial/\partial x_i}
\partial/\partial x_j=0$ for all $i,j=1,...,n$. We consider the
vector field $\xi$ on $\mathbb{R} ^n$ along $\mathbb{R} ^{n-1}$
defined by $\xi(x)=x_1 \partial/\partial x_n (x)$, where $x_1$ is
the first component of $x$. Now, to show that $\xi(\mathbb{R}
^{n-1})$ is totally geodesic in $T \mathbb{R} ^n$, it suffices to
check that (\ref{A19}) is verified. In fact,$\bar
\nabla_{\partial/\partial x_i} \bar \nabla_{\partial/\partial
x_j}\xi=\bar \nabla_{\partial/\partial x_i} \delta_{1j}
\partial/\partial x_n =0$. But $\bar \nabla_{\partial/\partial
x_1}\xi=\partial/\partial x_n$, and so $\xi$ is not parallel.
\end{remarks}

\section{The case of Lie groups with bi-invariant metrics}

Let us consider a connected Lie group $G^n$ equipped with a
bi-invariant metric $\bar g$, i.e. invariant by both left and
right translations. We shall generalize the results of Walczak~P.
\cite{Wk} on totally geodesic left invariant vector fields on
$G^n$ to left invariant vector fields along Lie subgroups.

Let $H^l$ be a Lie subgroup of $G^n$. The metric $g$ induced from
$\bar g$ on $H^l$ is a bi-invariant metric. If we denote by $\bar
\nabla$ and $\nabla$ the Levi-Civita connections on $G^n$ and
$H^l$ respectively, then we have $\bar \nabla_X
Y=\frac{1}{2}[X,Y]$, for all $X$,$Y$ of $\mathfrak{g}$, the Lie
algebra of $G^n$, and $\nabla_X Y=\frac{1}{2}[X,Y]$, for all
$X$,$Y$ of $\mathfrak{h}$, the Lie algebra of $H^l$.

\begin{lemma}\label{Subalg}
A connected complete submanifold $F^l$ of $G^n$ containing the
identity element $e$ of $G^n$, such that $T_e F^l$ is a subalgebra
of $\mathfrak{g}$, is totally geodesic if and only if $F^l$ is a
Lie subgroup $H^l$ of $G^n$.
\end{lemma}
\begin{proof}
If we denote by $\exp$ the exponential mapping $\exp:\mathfrak{g}
\to G^n$ of the Lie group $G^n$, and by $\exp_x:T_x G^n \to G^n$
the exponential map at a point $x$ of $G^n$ with respect to the
Levi-Civita connection of the metric $g$, then for all $x\in G^n$,
$\exp_x=\exp \circ (L_{x^{-1}})_*$, where $L_x$ is the left
translation of $G^n$ by $x$. Indeed, we show firstly that
$\exp_e=\exp$. Let $X\in \mathfrak{g}\equiv T_eG^n$ and
$\gamma(t)=\exp tX$. It suffices to check that $\gamma$ is a
geodesic. We have $\dot \gamma (t)=(L_{\gamma(t)})_*(\dot
\gamma(0))=(L_{\gamma(t)})_*(X)$, and thus $\bar \nabla_{\dot
\gamma (t)}\dot \gamma (t)=\bar
\nabla_{X(\gamma(t))}X(\gamma(t))$, where $X$ denotes also the
left invariant vector field on $G^n$ corresponding to $X$. Hence
$\bar \nabla_{\dot \gamma (t)}\dot \gamma
(t)=\frac{1}{2}[X,X](\gamma(t))=0$, and so $\exp_e=\exp$. Now, our
assertion follows from the fact that left translations are
isometries.

We consider a Lie subgroup $H^l$ of $G^n$ and $\mathfrak{h} =T_e
H^l$ its Lie algebra. If $X\in\mathfrak{h}$, then $\exp_e tX=\exp
tX \in H^l$, for all $t$ in a neighborhood of $0$, i.e. $H^l$
contains the geodesic starting from $e$ and with initial condition
$X$, and by the left translations, $H^l$ contains all geodesics
starting from points of $H^l$ with initial vectors tangent to
$H^l$ at these points. Thus $H^l$ is totally geodesic.

Conversely, suppose that $F^l$ is a connected complete
submanifold of $G^n$ such that $e\in F^l$ and $T_e
F^l=:\mathfrak{h}$ is a Lie subalgebra of $\mathfrak{g}$. Let
$H^l$ be the connected subgroup of $G^n$ with Lie algebra
$\mathfrak{h}$. $H^l$ is then a connected totally geodesic
submanifold of $G^n$ with $T_e H^l=T_e F^l$. Therefore $H^l=F^l$.
\end{proof}

\begin{proposition}\label{L-Invar}
A left invariant vector field on $G^n$ along a submanifold $F^l$
generates a totally geodesic submanifold of $TG^n$ if and only if
it is parallel along $F^l$ and $F^l$ is totally geodesic.
\end{proposition}
\begin{proof}
A left invariant vector field on $G^n$ is necessarily of constant
length, and we apply Theorem \ref{Th1}.
\end{proof}

\begin{corollary}
A left invariant vector field $\xi$ on $G^n$ along a Lie subgroup
$H^l$ is totally geodesic if and only if it is an element of the
centralizer of $\mathfrak{h}$ in $\mathfrak{g}$.
\end{corollary}
\begin{proof}
By Lemma \ref{Subalg}, $H^l$ is a totally geodesic submanifold in
$G^n$. Thus, by virtue of Proposition \ref{L-Invar}, $\xi$ is
totally geodesic if and only if $\xi$ is parallel along $H^l$.

Suppose that $\xi$ is totally geodesic. Then $\bar \nabla_X \xi
=0$, for all $X \in \mathfrak{h}$; i.e. $\xi$ is in the
centralizer of $\mathfrak{h}$ in $\mathfrak{g}$.

Conversely, if $\xi$ is in the centralizer of $\mathfrak{h}$ in
$\mathfrak{g}$, then $\bar \nabla_X \xi =0$, for all $X \in
\mathfrak{h}$. Let $x \in H^l$ and $z \in T_x H^l$. It suffices to
prove that $\bar \nabla_z \xi =0$. But $X:=(L_{x^{-1}})_*(z) \in
T_e H^l \equiv \mathfrak{h}$, and consequently $\bar \nabla_z \xi=
(\bar \nabla _X \xi)(x) =0$.
\end{proof}

\begin{corollary}
(a)\ There are no non-zero left invariant totally geodesic vector fields on a
semi-simple Lie subgroup of a Lie group with a bi-invariant Riemannian metric.

(b)\ Every left invariant vector field along a subgroup of an abeli\-an Lie group with
a bi-invariant Riemannian metric generates a totally geodesic submanifold of the
tangent bundle.

\end{corollary}

\begin{theorem}\label{CompCon}
Let $N^l$ be a connected complete totally geodesic embedded
submanifold of the tangent bundle of a connected Lie group $G^n$
equipped with a bi-invariant Riemannian metric such that
$H^l=\pi(N^l)$ is a Lie subgroup of $G^n$. Suppose that $N^l$ is
horizontal at a point $z$ of $T_e G^n$.

(a)\ If $z \in T_e H^l$, then $N^l$ is the image of $H^l$ by a left invariant vector
field on $H^l$ which belongs to the center of $\mathfrak{h}$. In particular, if $H^l$
is semi-simple, then $H^l$ is the only connected totally geodesic embedded submanifold
of $TG^n$ which is tangent to $H^l$ at $e$ and orthogonal to the fiber at a point of
$T_eG^n$.

(b)\ If $H^l$ is simple, then $N^l$ is the image of $H^l$ by a left invariant vector
field on $G^n$ along $H^l$ which belongs to the centralizer of $\mathfrak{h}$ in
$\mathfrak{g}$.
\end{theorem}
\begin{proof}
Using Proposition \ref{Transv}, there is a neighborhood $U$ of $e$
in $G^n$, a neighborhood $V$ of $z$ in $TG^n$ and a vector field
$Y$ on $M^n$ along $H^l \cap U$ such that $N^l \cap V=Y(H^l \cap
U), Y(e)=z$. We have $T_z N^l=T_z (N^l \cap V)=T_z Y(H^l \cap U)$.
Then each vector of $T_z N^l$ can be written as $X^h + (\bar
\nabla_X Y)^v$, for some $X\in \mathfrak{h}$. But $T_z N^l$ is a
subset of the horizontal subspace of $TTG^n$ at $z$, so at $e$ we
have $\bar \nabla_X Y=0$ for all $X\in \mathfrak{h}$. On the other
hand, since $N^l \cap V=Y(H^l \cap U)$ is totally geodesic, the
second assertion of Proposition \ref{Pr5} reduces at $e$ to the
identity $$ \bar \nabla_{X_1}\bar \nabla_{X_2}Y=\frac{1}{2}\bar
R(X_1,X_2)Y , \mbox{ for all vector fields } X_1,X_2 \mbox{ on }
\, H^l. $$ Then for all $W\in \mathfrak{g}=T_e G^n$, we have $$
\bar g(\bar \nabla_{X_1(e)}\bar \nabla_{X_2}Y,W)= \frac{1}{2} \bar
g(\bar R(X_1(e),X_2(e))Y(e),W). $$ If we extend $W$ to a vector
field $X_3$ along $H^l$, which is orthogonal to $\bar
\nabla_{X_2}Y$ in a neighborhood of $e$ in $H^l$, then we can
write $$ \bar g(\bar \nabla_{X_1(e)}\bar \nabla_{X_2}Y,W)=-\bar
g(\bar \nabla_{X_2(e)}Y,\bar \nabla_{X_1(e)}X_3)=0, $$ and
consequently, $ \bar g(\bar R(X_1(e),X_2(e))Y(e),W)=0, $ for all
$X_1(e),X_2(e)\in \mathfrak{h}=T_e H^l$ and $W\in \mathfrak{g}=T_e
G^n$. Therefore we have $$R(\cdot,\cdot)Y(e)=0, \mbox { when
applied to vectors in $T_e H^l$}. $$

Let us denote by $\xi$ the left invariant vector field on $G^n$
along $H^l$ such that $Y(e)=\xi(e)$. Then $\bar
R(\cdot,\cdot)\xi(e)=0$ when applied to vectors in $T_e H^l$, and
hence
\begin{equation}\label{center}
\bar R(\cdot,\cdot)\,\xi=0, \mbox{ when applied to elements of
$\mathfrak{h}$.}
\end{equation}

Consider now two cases.

(a)\ If $\xi(e)=z \in T_e H^l$, then $\xi \in \mathfrak{h}$, and we have, by virtue of
(\ref{center}), $\bar R(X,\xi)\xi=0$, for all $X \in \mathfrak{h}$. Thus
$|\,[\xi,X]\,|\,^2=4\bar g(\bar R(\xi,X)X,\xi)=0$ for all $X\in \mathfrak{h}$. It
follows that $\xi$ belongs to the center of $\mathfrak{h}$.

(b)\ If $H^l$ is simple, then $[\mathfrak{h},\mathfrak{h}]= \mathfrak{h}$. But $\bar
\nabla_{[X_1,X_2]} \xi= \frac12 [[X_1, X_2], \xi]= -2R(X_1,X_2)\xi=0$, for all $X_1$,
$X_2 \in \mathfrak{h}$, by virtue of (\ref{center}). Since
$[\mathfrak{h},\mathfrak{h}]= \mathfrak{h}$, we deduce easily that $\bar \nabla _X
\xi=0$, for all $X \in \mathfrak{h}$, or equivalently $[X,\xi]=0$, for all $X \in
\mathfrak{h}$. It follows that $\xi$ belongs to the centralizer of $\mathfrak{h}$ in
$\mathfrak{g}$.

In both cases, $\xi$ belongs to the centralizer of $\mathfrak{h}$ in $\mathfrak{g}$.
Hence, by Lemma \ref{Subalg}, $H^l$ is totally geodesic in $G^n$, and Proposition
\ref{L-Invar} implies then that $\xi(H^l)$ is a complete totally geodesic submanifold
of $TG^n$. Therefore $\xi(H^l)=N^l$, because $ \xi_* (T_e H^l)=T_z N^l$ and $N^l$ and
$H^l$ are connected.
\end{proof}

\begin{corollary}
Let $N^l$ be a connected complete horizontal totally geodesic
submanifold of the tangent bundle of a connected Lie group $G^n$
equipped with a bi-invariant Riemannian metric such that
$H^l=\pi(N^l)$ is a simply connected submanifold of $G^n$
containing the identity element. Suppose that $\mathfrak{h}:=
\pi_*(T_z N^l)$ is a Lie subalgebra of $\mathfrak{g}$ for a point
$z$ of $T_e G^n\cap N^l$. If $Z \in T_e H^l$ (resp. $\mathfrak{h}$
is simple), then $H^l$ is a Lie subgroup of $G^n$ and $N^l$ is the
image of $H^l$ by a left invariant vector field on $H^l$ (resp. on
$G^n$ along $H^l$) which belongs to the center of $\mathfrak{h}$
(resp. centralizer of $\mathfrak{h}$ in $\mathfrak{g}$).
\end{corollary}
\begin{proof}
By Theorem \ref{Hor-Sub}, $H^l$ is complete and totally geodesic.
It follows from Lemma \ref{Subalg} that $H^l$ is a Lie subgroup of
$G^n$. Now, our corollary follows from Theorem \ref{CompCon}.
\end{proof}

\vspace{1ex}

\noindent
Mohamed Tahar Kadaoui ABBASSI,\\
D\'epartement des Math\'ematiques, \\
Facult\'e des sciences Dhar El Mahraz,\\
Universit\'e Sidi Mohamed Ben Abdallah,\\
B.P. 1796, Fes-Atlas,\\
Fes, Morocco\\
e.mail: mtk{\_}abbassi@Yahoo.fr\\
\vspace{1ex}

\noindent
Alexander YAMPOLSKY,\\
Department of Geometry,\\
Faculty of Mechanics and Mathematics,\\
Kharkiv National University,\\
Svobody Sq. 4,\\
61077, Kharkiv,\\
Ukraine.\\
e-mail: yamp@univer.kharkov.ua

\begin{thebibliography}{30}
\bibitem{BX-V1}
Boeckx~E., Vanhecke~L. Harmonic and minimal radial vector fields,
Acta Math. Hungar. 90(2001), 317-331.

\bibitem{BX-V2}
Boeckx~E., Vanhecke~L. Harmonic and minimal vector fields on
tangent  and unit tangent bundles, Differential Geom. Appl.13
(2000), 77-93.


\bibitem{Dmb} Dombrowski P. {On the geometry of tangent bundle}, J. Reine
Angew. Math., 210 (1962), N~1-2, 73-88.

\bibitem{GMLF}Gil-Medrano~O., Llinares-Fuster~E. Second variation
of volume and energy of vector fields. Stability of Hopf vector
fields, Math. Ann. 320 (2001), 531-545.

\bibitem{GM} Gil-Medrano~O., Llinares-Fuster~E. Minimal unit vector
fields, T\^ohoku Math. J., 54(2002), 71-84.

\bibitem{G-Z}
Gluck~H., Ziller~W. On the volume of a unit vector field on the
three-sphere, Comm. Math. Helv. 61 (1986), 177-192.

\bibitem{GD-V1}
Gonz\'alez-D\'avila J.C., Vanhecke L. Examples of minimal unit
vector fields, Ann. Global Anal. Geom. 18 (2000), 385-404.

\bibitem{KoNz}
Kobayashi S., Nomizu K. \emph{Foundations of differential
geometry}, Interscience Publ., Vol.1 (1967) and Vol.1,2 (1969).

\bibitem{Kow}
Kowalski O. {Curvature of the induced Riemannian metric on the
tangent bundle of a Riemannian manifold}, J. Reine Angew. Math.
250 (1971), 124-129.

\bibitem{Liu}
Liu M.-S. { Affine maps of tangent bundles with Sasaki metric},
Tensor, N.S., 28 (1974), 34-42.

\bibitem{Rov} Rovenskii V. {\em Foliations on Riemannian Manifolds
and Submanifolds.} Birkh\"auser, 1997.

\bibitem{Sk}
Sasaki~S. {On the differential geometry of tangent bundles of
Riemannian manifolds}, T\^ohoku Math. J. 10 (1958), 338-354.

\bibitem{St}
Sato K.{ Geodesics on the tangent bundles over space forms}, Tensor,
32 (1978), 5-10.

\bibitem{Wk}
Walczak P. {On totally geodesic submanifolds of tangent bundles
with Sasaki metric}, Bull. Acad. Pol. Sci., ser. Sci. Math., 28
(1980), N~3-4, 161-165.

\bibitem{Wk-Bt} Walczak P. On the energy of unit vector fields with
isolated singularities, Ann. Pol. Math., LXXII.3 (2000), 269-274.

\bibitem{Wr} Warner, F.W. {\em Foundations of differentiable
manifolds and Lie groups}, Academic Press, New York, 1971.

\bibitem{Wd} Wood C.M. On the energy of a unit vector field,
Geom. Dedicata 64 (1997), 319-330.

\bibitem{Ym1} Yampolsky A. On the mean curvature  of a unit vector
field,{ Math. Publ. Debrecen, 60/1-2 (2002), 131 -- 155.}

\bibitem{Ym2} Yampolsky A. On the intrinsic geometry of a unit
vector field, Comment. Math. Univ. Carolinae 43, 2 (2002),
299-317.

\bibitem{Ym} Yampolsky A. Totally geodesic property of the Hopf
vector field, Acta Math. Hungarica, Acta   Math. Hungar. 101, 1-2 (2003), 73-92.   
\end{thebibliography}
\end{document}